\documentclass[12pt]{amsart}

\usepackage{amsmath,amssymb,amsfonts,fullpage,amssymb,amsthm}

\usepackage{verbatim}
\usepackage[usenames]{color}
\usepackage{hyperref}
\usepackage{url}
\usepackage{tikz,tikz-qtree,ifthen,cancel}
\usepackage{array,tikz-qtree,ifthen,cancel}
\usepackage{graphicx}
\usepackage{adjustbox}
\usepackage{amsthm,graphicx,tikz,appendix,tikz-qtree,ifthen,cancel}
\usetikzlibrary{calc,shapes,patterns,positioning}

\theoremstyle{remark}

\theoremstyle{definition}

\theoremstyle{remark}

\newcommand{\om}{\omega}

\newcommand{\Lra}{\Longrightarrow}

\newcommand{\Lauchli}{L{\"{a}}uchli}

\newcommand{\noprint}[1]{\relax}

\title[Problems on Reverse Mathematics of Ramsey Theory]{A List of Problems on the Reverse Mathematics of Ramsey Theory on the Rado Graph and\\ on Infinite,  Finitely Branching Trees}

\author{Natasha Dobrinen}
\address{Department of Mathematics\\
 University of Denver \\
C.M.\ Knudson Hall, Room 300\\
2390 S.\ York St.\\ Denver, CO \ 80208 U.S.A.}
\email{natasha.dobrinen@du.edu}
  \urladdr{\url{http://web.cs.du.edu/~ndobrine}}

\begin{document}

\maketitle

\begin{abstract}
This list presents 
 problems in the Reverse Mathematics of  infinitary Ramsey theory which I find interesting but do not personally have the techniques to solve. 
The intent is to enlist the help of  those working in Reverse Mathematics to take on such problems, and the myriad of related questions one can infer from them.
A short bit of background and starting references are provided.
\end{abstract}

\section{Motivation}

In \cite{Csima/Mileti09}, Csima and Mileti proved that the reverse mathematical strength of the rainbow Ramsey theorem for pairs of natural numbers is strictly weaker than the Ramsey theorem for pairs of natural numbers, over RCA$_0$.  We wonder how much will remain the same or be different when a graph structure is added, and what the relationships are between those and the relevant Ramsey theorems on infinite trees?

The {\em Rado graph} $\mathcal{R}$  is the  homogeneous graph  on countably many vertices which is universal for all graphs on countably many vertices.  It is also called the random graph, as it can be constructed by flipping a coin to decide when to add an edge between two vertices.
Pouzet and Sauer proved  that the Rado graph does not have the Ramsey property.  
In particular, they showed that there is an edge coloring with two colors such that no subgraph which is isomorphic to the original (and thus again a Rado graph) has all edges of the same color.  
See \cite{Pouzet/Sauer96}.

However, the Rado graph does have {\em finite big Ramsey degrees}, proved in \cite{Sauer06} and \cite{Laflamme/Sauer/Vuksanovic06}.
This means that given any finite graph $G$, there is a number bound
$T(G)$ such that for any coloring of all copies of $G$ in $\mathcal{R}$ into finitely many colors, there is a subgraph $\mathcal{R}'$ of $\mathcal{R}$, which is again a Rado graph, in which all copies of $G$ take no more than $T(G)$ colors. 
The unavoidable color classes are called {\em canonical partitions}, and are described in terms of tree structures as subtrees of $2^{<\om}$.
Their proof uses  a Ramsey theorem of Milliken for strong trees \cite{Milliken79} and \cite{Milliken81}, which has at its core the Halpern-\Lauchli\ Theorem \cite{Halpern/Lauchli66}.
More background on these theorems can be found in Chapters 3 and 6 of Todorcevic's book \cite{TodorcevicBK10}, and also in my expository paper \cite{DobrinenRIMS17} available on my website.

In \cite{Dobrinen/Laflamme/Sauer16}, Dobrinen, Laflamme, and Sauer proved that  the Rado graph has the rainbow Ramsey property.
This means that fixing any finite graph $G$ and $k\ge 2$, given any coloring of all copies of $G$ in $\mathcal{R}$ into $\om$ many colors such that  each color appears no more than $k$ times,
there is a subgraph $\mathcal{R}'$ of $\mathcal{R}$ which is again a Rado graph, and such that no two copies of $G$ in $\mathcal{R}'$ have the same color.  Some reverse mathematics questions appear in \cite{Dobrinen/Laflamme/Sauer16}, and we augment that list here.

The following implications are known as ZFC, or even just ZF, results: with (A) $\Lra$ (B) $\Lra$ (C) $\Lra$ (D).

\begin{enumerate}
\item[(A)]
The Halpern-\Lauchli\ Theorem for any finite number of trees.   This is a ZF theorem.
\item[(B)]
Milliken's Ramsey theorem for strong trees (on one finitely branching infinite tree with no terminal nodes).  It is proved by induction using the Halpern-\Lauchli\ Theorem.
\item[(C)]
Canonical partition theorem on colorings of  a fixed finite graph inside the Rado graph. 
(This implies that the Rado graph has finite big Ramsey degrees.)
\item[(D)]
The Rado graph has the rainbow Ramsey property, for $k$-bounded colorings, for each $k\ge 2$.

\end{enumerate}

\section{Problems in the Reverse Mathematics of the Rado Graph and Related Ramsey Theorems on Infinite, Finitely Branching Trees}

\begin{enumerate}
\item
Are any of (A) - (D) equivalent over RCA$_0$?

\item
Is the rainbow Ramsey theorem for edges in the Rado graph strictly stronger than the rainbow Ramsey theorem for pairs on the natural numbers?

\item
What about $k$-bounded colorings of graphs with three vertices in the Rado graph versus $k$-bounded colorings of triples of natural numbers?
That is, assuming the rainbow Ramsey theorem for $k$-bounded colorings of triples of natural numbers, does the rainbow Ramsey theorem for $k$-bounded colorings of copies of $G$ in the Rado graph follow over RCA$_0$, where $G$ is any graph on three vertices?

\item 
Is the reverse math strength of the (strong tree version) Halpern-\Lauchli\ Theorem for $d$ trees strictly weaker than for $d+1$-trees?

\item 
Is the reverse math strength of the strong tree version of Halpern-\Lauchli\ for $d$ trees strictly stronger than the somewhere dense matrix form for $d$ trees?
(This would contrast with results in \cite{Dobrinen/Hathaway16} (see also \cite{Dobrinen/Hathaway17Addendum}), where we showed that these two forms are equivalent for strongly inaccessible cardinals, over ZFC.)

\item
Does ``The Halpern-\Lauchli\ Theorem holds for any finite number of trees" have  the same strength as ``Milliken's Theorem for colorings of finite  strong trees holds"?

\item 
Does the strength of the Milliken Theorem depend  on the height of the finite trees being colored?
\item
Does the strength of either the Halpern-\Lauchli\ Theorem or the Milliken Theorem depend on the splitting properties of the trees? 
For instance, are two-branching trees weaker than three-branching trees,..., weaker than finitely but unboundedly branching trees?

\item  
What is the Reverse Mathematics strength of the statement  (D)?

\item 
Does the rainbow Ramsey property for $k$-bounded colorings of the Rado graph have the same reverse math strength as  for $(k+1)$-bounded colorings?

\item
Does the rainbow Ramsey theorem for the Rado graph have the same reverse math strength as the canonical partition theorem of Sauer?

\item
How do the strengths of the rainbow Ramsey theorem for the Rado graph compare for colorings of different finite graphs?
Are they all the same, or does the number of vertices in the finite graph being colored affect the strength?

\end{enumerate}

\bibliographystyle{amsplain}
\bibliography{references}

\providecommand{\bysame}{\leavevmode\hbox to3em{\hrulefill}\thinspace}
\providecommand{\MR}{\relax\ifhmode\unskip\space\fi MR }
\providecommand{\MRhref}[2]{%
  \href{http://www.ams.org/mathscinet-getitem?mr=#1}{#2}
}
\providecommand{\href}[2]{#2}
\begin{thebibliography}{10}

\bibitem{Csima/Mileti09}
Barbara Csima and Joseph Mileti, \emph{The strength of the rainbow {R}amsey
  theorem}, Journal of Symbolic Logic \textbf{74} (2009), no.~4, 1310--1324.

\bibitem{DobrinenRIMS17}
Natasha Dobrinen, \emph{Forcing in {R}amsey theory}, Proceedings of 2016 RIMS
  Symposium on Infinite Combinatorics and Forcing Theory (2017), 17--33.

\bibitem{Dobrinen/Hathaway16}
Natasha Dobrinen and Daniel Hathaway, \emph{The {H}alpern-{L}{\"{a}}uchli
  {T}heorem at a measurable cardinal}, Journal of Symbolic Logic \textbf{82}
  (2017), no.~4, 1560--1575.

\bibitem{Dobrinen/Hathaway17Addendum}
\bysame, \emph{Addendum to the {H}alpern-{L}{\"{a}}uchli {T}heorem at a
  measurable cardinal}, Journal of Symbolic Logic (2018), To appear.

\bibitem{Dobrinen/Laflamme/Sauer16}
Natasha Dobrinen, Claude Laflamme, and Norbert Sauer, \emph{Rainbow {R}amsey
  simple structures}, Discrete Mathematics \textbf{339} (2016), no.~11,
  2848--2855.

\bibitem{Halpern/Lauchli66}
J.~D. Halpern and H~L{\"{a}}uchli, \emph{A partition theorem}, Transactions of
  the American Mathematical Society \textbf{124} (1966), 360--367.

\bibitem{Laflamme/Sauer/Vuksanovic06}
Claude Laflamme, Norbert Sauer, and Vojkan Vuksanovic, \emph{Canonical
  partitions of universal structures}, Combinatorica \textbf{26} (2006), no.~2,
  183--205.

\bibitem{Milliken79}
Keith~R. Milliken, \emph{A {R}amsey theorem for trees}, Journal of
  Combinatorial Theory, Series A \textbf{26} (1979), 215--237.

\bibitem{Milliken81}
\bysame, \emph{A partition theorem for the infinite subtrees of a tree},
  Transactions of the American Mathematical Society \textbf{263} (1981), no.~1,
  137--148.

\bibitem{Pouzet/Sauer96}
Maurice Pouzet and Norbert Sauer, \emph{Edge partitions of the {R}ado graph},
  Combinatorica \textbf{16} (1996), no.~4, 505--520.

\bibitem{Sauer06}
Norbert Sauer, \emph{Coloring subgraphs of the {R}ado graph}, Combinatorica
  \textbf{26} (2006), no.~2, 231--253.

\bibitem{TodorcevicBK10}
Stevo Todorcevic, \emph{Introduction to {R}amsey {S}paces}, Princeton
  University Press, 2010.

\end{thebibliography}

\end{document}